\documentclass[10pt]{amsart}
\usepackage[utf8]{inputenc}

\usepackage{amsmath}
\usepackage{amssymb}
\usepackage{amsthm}
\usepackage{amsfonts}
\usepackage{mathrsfs}
\usepackage{commath}
\usepackage{tikz}
\usepackage{cancel}
\usepackage{graphicx}
\usepackage{enumerate}

\usepackage{ifthen}

\usepackage{listings}

\usepackage{hyperref,url,pst-node}
\hypersetup{citecolor=red, linkcolor=blue, colorlinks=true}

\usepackage{mathtools}

\newcommand{\CANTOR}[3]{
\ifnum 0 < \numexpr#1 {
    \CANTOR{\numexpr#1-1}{#2*0.333333}{#3}
    \CANTOR{\numexpr#1-1}{#2*0.333333}{#3+0.666666*#2}}
\else{
    \draw[black,line width=0.5cm] (#3,0)--(#3+#2,0);}
\fi{}
}

\theoremstyle{plain}

\newtheorem{theorem}{Theorem}[section]

\newtheorem{proposition}{Proposition}
\newtheorem*{proposition*}{Proposition}

\newtheorem*{conjecture*}{Conjecture}
\theoremstyle{definition}

\theoremstyle{remark}

\title{The Entropy and Hausdorff Dimension of self-similar sets}
\author{James Evans}
\date{September 2020}

\begin{document}

\maketitle

\begin{abstract}
    Given a $k$-self similar set $X\subset [0,1]^{d}$ we calculate both its Hausdorff dimension and its entropy, and show that these two quantities are in fact equal. This affirmatively resolves a conjecture of Adamczewski and Bell \cite{BellAutomaticFractals}.
\end{abstract}

\section{Introduction}


A $k$-automatic set is a subset $\mathcal{N}$ of $\mathbb{N}$ whose elements, $n$, are recognisable by a deterministic finite automaton taking the base $k$ expansion of $n$ as input. Equivalently (by Eilenberg, \cite{EilenbergAutomata}), $\mathcal{N}$ is $k$-auomatic if the collection of all sets of the form

$$
K^{k}_{a,b} = \{n \geqslant 0\; | \; k^{a}n+b \in \mathcal{N}\}, \quad a \geqslant 0, \quad 0 \leqslant b < k^{a}
$$
is finite\footnote{See Allouche and Shallit \textit{Automatic sequences} \cite{AutomaticSequences} for all the details.}. These sequences, their generalisations and families of related objects (which include the $k$-regular sequences and the $k$-Mahler functions), are prominent in many different areas of mathematics.

One important common feature of these related families of objects is that they exhibit independence results. The prototypical example is the result proven by Cobham \cite{CobhamTheorem}, that a set is $k$-automatic and $l$-automatic for multiplicatively independent $k$ and $l$ ($k^{a} \neq l^{b}$ for all $a,b \in \mathbb{N}$) if and only if it is eventually periodic. Later Bell \cite{BellRegularSequences} proved the analogous result for $k$-regular sequences, and Adamczewski and Bell proved the one for Mahler functions \cite{AdamczewskiBellMahlerFunctions}.

As part of their ongoing effort to generalise and study Cobham-type independence results, Adamczewski and Bell \cite{BellAutomaticFractals} recently defined yet another generalisation of $k$-automatic subsets of the integers: the $k$-self similar subsets of $[0,1]^{d}$. These are defined as follows. Let $X \subset [0,1]^{d}$ be compact, and fix an integer $k \geqslant 2$. Let $a \geqslant 0$, $\bold{b}=(b_{1},\ldots,b_{d})\in \mathbb{Z}^{d}$ such that $0 \leqslant b_{i} < k^{a}$ for each $i$, and consider the subsets $\mathcal{K}_{a,\bold{b}}^{k}$ of $X$ defined by
$$
\mathcal{K}_{a,\bold{b}}^{k}:=\{(k^{a}x_{1}-b_{1},\ldots,k^{a}x_{d}-b_{d}) \; | \; (x_{1},\ldots,x_{d}) \in X \cap \prod_{i=1}^{d}[\frac{b_{i}}{k^{a}},\frac{b_{i}+1}{k^{a}}]\}.
$$
That is, takes the $\frac{1}{k^{a}}$-hyper-cube $\prod_{i=1}^{d}[\frac{b_{i}}{k^{a}},\frac{b_{i}+1}{k^{a}}]$, intersect it with $X$ and scale this intersection back up to $[0,1]^{d}$. The $k$-kernel of $X$ is the collection of all such sets

$$
\mathcal{K}_{X}^{k}=\{\mathcal{K}_{a,\bold{b}}^{k} \; : 0 \leqslant a, \; 0 \leqslant b_{i} < k^{a}\; \}.
$$
If $\mathcal{K}_{X}^{k}$ is finite then we say that $X$ is a \emph{$k$-self-similar set}. This definition is designed to closely mirror that for $k$-automatic subsets of the integers.

Adamczewski and Bell \cite{BellAutomaticFractals} (for $d=1$), and Chan and Hare \cite{MultiDimensionalFractalCobham} (for all $d$) have proven a Cobham-type result for these objects: $X$ is both $k$- and $l$-self-similar for multiplicatively independent $k$ and $l$ if and only if it is a union of finitely many rational polyhedra. In this case, $X$ will has Hausdorff dimension $0$ or $1$.

For this (and other reasons) it is interesting to calculate the Hausdorff dimension of the set $X$. The paper \cite{BellAutomaticFractals} includes an additional conjecture by the authors in this direction, connecting the Hausdorff dimension of $X$ to its so called \textit{entropy}. In fact, they only make the conjecture for the subclass of $k$-automatic sets, and not all $k$-self-similar sets. However, with a generalisation of their definition of entropy we prove the result for all self-similar sets.
\begin{theorem} \label{Theorem}
For every $k$-self-similar set $X \subset [0,1]^{d}$, we have $\text{\rm{dim}}_{\mathcal{H}}(X)=h(X)$,
where $\text{\rm{dim}}_{\mathcal{H}}(X)$ is the Hausdorff dimension of the set, and $h(X)$ is its entropy.
\end{theorem}
We shall have more to say about this result and its connections in the conclusion.

\section{Proof of Theorem \ref{Theorem}}

In this section we prove Theorem \ref{Theorem}. For convenience, fix the notation so that $X$ always refers to a self-similar set, $k$ its index of self-similarity and $d$ the dimension. We first calculate the Hausdorff dimension of $X$ by using a method of Mauldin and Williams \cite{MauldinWilliams}. Subsequently, we shall define and then calculate the entropy of $X$. A comparison of the two resulting formulas yields Theorem \ref{Theorem}.

In order to calculate the Hausdorff dimension\footnote{ 
For the required theory of fractals and Hausdorff dimension, see Falconer \cite{FalconerFractalGeometry}.} of $X$, we utilise the work of Mauldin and William \cite{MauldinWilliams} and an example by Marion \cite{MarionExample} contained therein, and realise $X$ using a geometric graph-directed construction\footnote{Recently Charlier, Leroy and Rigo \cite{AnalogueOfCobhamForGDIFS} linked and generalised a variety of Cobham-like results using the closely related concept of Graph Directed Iterated Function Systems. There are interesting connections between the work of that paper and the current one.}. These are a vast generalisation of the iterated function system concept from fractal geometry. A geometric graph-directed construction in $\mathbb{R}^{d}$ consists of the following.
\begin{itemize}
    \item A finite collection $J_{0},\ldots,J_{n-1}$ of compact subsets of $\mathbb{R}^{d}$ with non-empty and non-intersecting interiors. 
    \item A directed graph $G$ which has similarity maps $T_{ij}$, with similarity ratios $t_{ij}$, associated to each edge $ij$. For each $i$
    \begin{itemize}
        \item there is some edge $ij\in E(G)$ out of vertex $i$.
        \item the $T_{ij}(J_{j})$, $ij \in E(G)$ are contained in $J_{i}$ and have disjoint interiors.
        \item If $[i_{0},i_{2},\ldots,i_{q}=i_{0}]$ is a directed cycle in $G$, then $\prod_{k=0}^{q-1}t_{i_{k}}t_{i_{k+1}}<1$.
    \end{itemize}
\end{itemize}
Let $A_{G}$ be the weighted adjacency matrix of the graph $G$, that is $A_{G}=\{t_{ij}\}_{i,j=0}^{n-1}$ (with $t_{ij}=0$ if there is no edge $ij$). Let $A_{G,\beta}=\{t_{ij}^{\beta}\}_{i,j=0}^{n-1}$ and $\rho(A_{G,\beta})$ be the spectral radius of $A_{G,\beta}$. For $p \geqslant 1$ denote the set of length $p$ paths in $G$ by $G(p)=\{\sigma=(\sigma(0),\ldots,\sigma(p))\in \{0,\ldots,n-1\}^{p+1} \; | \; \sigma(j)\sigma(j+1) \in E(G), j=0,\ldots,p-1\}$. 

\begin{theorem}[Mauldin \& Williams \cite{MauldinWilliams}] \label{GGDC}
For each geometric graph-directed construction there is a unique list $K_{0},\ldots,K_{n-1}$ of compact subsets of $\mathbb{R}^{d}$ such that $K_{i}=\bigcup_{ij \in E(G)} T_{ij}(K_{j})$. The construction object $K= \bigcup_{i=0}^{n-1}K_{i}$ is equal to
$$K=\bigcap_{p=1}^{\infty} \bigcup_{\sigma \in G(p)} T_{\sigma}(J_{\sigma(p)}),$$
where $T_{\sigma}=T_{\sigma(0)\sigma(1)}\circ \cdots \circ T_{\sigma(p-1)\sigma(p)}$. It has Hausdorff dimension $\alpha$, where $\alpha$ is the unique value such that $\rho(A_{G,\alpha})=1$.
\end{theorem}

Let $\mathcal{K}_{X}^{k}=\{X_{0},\ldots,X_{n-1}\}$ be the $k$-kernel of $X$. Without loss of generality, $X=X_{0}$. We shall define a GGDC which realises $\textbf{X}=\bigcup_{i=0}^{n-1} (X_{i}+\{(2i,0,\ldots,0)\})$ as its construction object. This exact form for $\textbf{X}$ is not crucial---we only need it to contain copies of each $X_{j}$ which do not overlap. The salient features of our set up, which can be proven from the definition of the $k$-kernel, are
\begin{enumerate}
    \item The $X_{j}$ are compact subsets of $[0,1]^{d}$.
    \item Each $X_{j}$ is the union of a collection of $\frac{1}{k}$ scaled copies of the other $X_{i}$. These copies are the intersections of $\frac{1}{k}$ grid hyper-cubes $\prod_{l=1}^{d}[\frac{b_{l}}{k},\frac{b_{l}+1}{k}]$ with $X_{j}$. \label{ref}
\end{enumerate}
Now let $A_{ij}$ be the number of $\frac{1}{k}$ scaled copies of $X_{i}$ contained in $X_{j}$. The vertex set of the graph $G$ is given by 
$$
V(G)=\{(h,i,j)|1\leqslant  h \leqslant A_{ij}, \; 0 \leqslant i,j \leqslant n-1\}.
$$
We imagine each vertex $(h,i,j)$ as corresponding to $X_{hij}$, the $h$'th copy of $X_{i}$ contained inside the copy of $X_{j}$ in $\textbf{X}$. Two vertices $(h,i,j)$, $(r,s,t)$ are connected by an edge if $i=t$. That is, we connect each $\frac{1}{k}$ scaled copy of $X_{i}$ in $\textbf{X}$ (each $X_{hij}$) to every piece contained within the original $X_{i}$ (every $X_{rsi}$). By \eqref{ref} each $X_{hij}$ is the intersection of $X_{j}$ with some closed $\frac{1}{k}$ grid hyper-cube. The corresponding set $J_{hij}$ associated to the vertex $(h,i,j)$ is this hyper-cube. The similarity map $T_{(hij),(rsi)}$ associated to an edge is the obvious similarity mapping from the big hyper-cube containing $X_{i}$ onto $J_{hij}$. This means that $T_{(hij),(rsi)}$ is independent of $r$ and $s$, and has scaling ratio $\frac{1}{k}$.

One can check that all of the conditions required of a geometric graph-directed construction are obeyed by this setup. The only task remaining is to show that the construction object $K=\bigcap_{p=1}^{\infty} ( \cup_{\sigma \in G(p)} T_{\sigma}(J_{\sigma(p)}) )$ actually equals $\textbf{X}$. The key point is to note that the $p$'th level of the construction object, $\cup_{\sigma \in G(p)} T_{\sigma}(J_{\sigma(p)})$ is the union of all $\frac{1}{k^{p+1}}$ grid hyper-cubes which are contained in some $J_{hij}$ and intersect $\textbf{X}$. Because $\textbf{X}$ is closed, any point not in $x$ is contained in some sufficiently small hyper-cube which does not intersect $X$ and this point is not in $K$. However, if $x$ is in $\textbf{X}$ then every hyper-cube containing it does intersect $\textbf{X}$, and $x$ is in $K$. Thus $K=\textbf{X}$.

The weighted adjacency matrix of the graph $G$ is $\frac{1}{k}\tilde{A}$, where $\tilde{A}$ is its unweighted adjacency matrix. Theorem \ref{GGDC} hence guarantees that $\textbf{X}$ has dimension $\log_{k}(\rho(\tilde{A}))$, $\rho(\tilde{A})$ the spectral radius of $\tilde{A}$. However, contained in \cite{MauldinWilliams} is a proof that $\tilde{A}$ and the matrix $A=\{A_{ij}\}_{i,j=0}^{n-1}$ have the same spectral radius. Thus $\text{dim}_{\mathcal{H}}(\textbf{X})=\log_{k}(\rho(A))$.

But $X=X_{0}$ is only a small piece of $\textbf{X}$, so in general its Hausdorff dimension will be less than or equal to that of $\textbf{X}$\footnote{In general not every $X_{i}$ will have the same dimension. Whether they do or not is related to whether the graph $G$ is strongly connected, which guarantees that every $X_{i}$ will contain some scaled copy of every $X_{j}$.}. However, we can show that $\text{ dim}_{\mathcal{H}}(\textbf{X})=\text{ dim}_{\mathcal{H}}(\textbf{X}_{0})$. Now, $\textbf{X}$ is the union of disjoint translated copies of the $X_{i}$ so $\text{dim}_{\mathcal{H}}(\textbf{X})=\max \{\text{dim}_{\mathcal{H}}(X_{i})\}$. But by definition, for each $i$, $X$ must contain (or is equal to) a $\frac{1}{k^{p}}$ scaled copy of $X_{i}$ for some $p$. Hence by monotonicity of dimension, $\text{dim}_{\mathcal{H}}(\textbf{X}) \geqslant \text{dim}_{\mathcal{H}}(X) \geqslant \max \{\text{dim}_{\mathcal{H}}(X_{i})\} = \text{dim}_{\mathcal{H}}(\textbf{X})$. 

The preceding considerations have thus proven the following proposition.

\begin{proposition}\label{1}
Let $X$ be a $k$-self-similar set and $A$ the matrix $A=\{A_{ij}\}_{i,j=0}^{n-1}$, $A_{ij}$ is the number of $\frac{1}{k}$ scaled copies of $X_{i}$ contained within $X_{j}$. Then
$$\text{dim}_{\mathcal{H}}(X)=\log_{k}(\rho(A)).$$
\end{proposition}

\bigskip

Now we discuss the entropy $h(X)$ of our self-similar set $X$. We begin with a very brief summary of some terminology from the study of combinatorics on words. 

Take any finite set $\mathcal{A}$. A word, $w$, over $\mathcal{A}$ is a formal concatenation of symbols from $\mathcal{A}$, $w=a_{1}a_{2}\cdots a_{n}$. The set of all words of finite length over $\mathcal{A}$ is denoted by $\mathcal{A}^{*}$. A language, $\mathcal{L}$, over $\mathcal{A}$ is a subset of $\mathcal{A}^{*}$. A language $\mathcal{L}$ is factorial if given any word $w$ in $\mathcal{L}$, every subword of $w$ is also in $\mathcal{L}$. Let $\mathcal{P}(\mathcal{L},n)$ denote the number of length $n$ words in $\mathcal{L}$. If $\mathcal{L}$ is factorial, then $\mathcal{P}(\mathcal{L},n+m)\leqslant \mathcal{P}(\mathcal{L},n)\mathcal{P}(\mathcal{L},m)$. The entropy of $\mathcal{L}$ (which exists if $\mathcal{L}$ is factorial) is defined to be
$$
h(\mathcal{L})=\lim_{n \to \infty}\frac{1}{n}\log(\mathcal{P}(\mathcal{L},n)).
$$

To calculate the entropy, $h(X)$ of our set $X$, we first associate to it a factorial language, $\mathcal{L}_{X}$. Then $h(X)$ is the entropy of $\mathcal{L}_{X}$, however we modify things slightly and instead use the base $k$-logarithm in the definition.
$$
h(X)=\lim_{n \to \infty}\frac{1}{n}\log_{k}(\mathcal{P}(\mathcal{L}_{X},n)).
$$
One reason for this modification is that it makes Theorem \ref{Theorem} true. If a set $X$ is $k$-self-similar, we can also consider it as $k^{a}$ self-similar for $a \geqslant 2$. However we want a definition of entropy which is independent of $k$, since we want it to be equal to the Hausdorff dimension which depends only on $X$. Taking the index of self-similarity as the base of the logarithm helps to guarantee this. 

To define $\mathcal{L}_{X}$, let $\mathcal{K}_{X}^{k}=\{X_{0},\ldots,X_{n-1}\}$ be the $k$-kernel of $X$ with $X=X_{0}$. The basic alphabet of $\mathcal{L}_{X}$ is the set of $d$-tuples $\varSigma_{k}^{d}=\{(b_{1},\ldots,b_{d}) \; | \; b_{i}=0,\ldots,k-1\}$. Now suppose that $x=(x_{1},\ldots,x_{d}) \in X$ and that each $x_{i}$ has a $k$-adic representation $x_{i}=(b_{i}^{(1)},b_{i}^{(2)},b_{i}^{(3)},\ldots)_{k}$. We pad the $k$-adic representation with infinitely many trailing zeroes if necessary. Then take the infinite word on the alphabet $\varSigma_{k}^{d}$
$$
\textbf{w}(x)=\begin{pmatrix}b_{1}^{(1)}\\b_{2}^{(1)}\\ \vdots \\b_{d}^{(1)} \end{pmatrix}\begin{pmatrix}b_{1}^{(2)}\\b_{2}^{(2)}\\ \vdots \\ b_{d}^{(2)} \end{pmatrix} \begin{pmatrix}b_{1}^{(3)} \\ b_{2}^{(3)} \\ \vdots \\ b_{d}^{(3)} \end{pmatrix}\cdots
$$
and add all finite subwords of $\textbf{w}(x)$, that is all subwords of the following form
$$
\textbf{w}=\begin{pmatrix}b_{1}^{(m+1)}\\b_{2}^{(m+1)}\\ \vdots \\ b_{d}^{(m+1)} \end{pmatrix} \begin{pmatrix}b_{1}^{(m+2)}\\b_{2}^{(m+2)}\\ \vdots \\ b_{d}^{(m+2)} \end{pmatrix}\cdots \begin{pmatrix}b_{1}^{(m+p)}\\b_{2}^{(m+p)}\\ \vdots \\ b_{d}^{(m+p)} \end{pmatrix}
$$
to $\mathcal{L}_{X}$ (this makes $\mathcal{L}_{X}$ a factorial language). We do this for all $x$ in $X$, and all $k$-adic representations. That is, if some $x^{(i)}$ is itself a $k$-adic rational we \textbf{do} include both of its $k$-adic representations. 

This definition closely follows that of Adamcewski and Bell \cite{BellAutomaticFractals}. The main change is that instead of looking at the base $k$ expansions of $k$-adic rationals in $X$, we look at the expansions of all elements of $X$. In general a $k$-self-similar set need not contain any $k$-adic rationals. If $X$ is a $k$-automatic set, this gives a slightly different language than \cite{BellAutomaticFractals} because we are including the multiple representations of the $k$-adic rationals. However, this will not affect the entropy, because including both representations will at worst increase $\mathcal{P}(\mathcal{L}_{X},n)$ by a constant multiple, which will not affect the asymptotics growth of the logarithm.

It can quite easily be shown, using the definition of the $k$-kernel of $X$, that we have the following equivalent characterisation of $\mathcal{L}_{X}$. Take any $d$-tuple of $k$-adic rational numbers which has $p$ $k$-adic digits, $((0.b_{1}^{(1)}\cdots b_{1}^{(p)})_{k},\ldots, (0.b_{d}^{(1)}\cdots b_{d}^{(p)})_{k}))$. Pad them with trailing zeroes if needed so they all have $p$ digits. The word
$$
\textbf{w}=\begin{pmatrix}b_{1}^{(1)}\\b_{2}^{(1)}\\ \vdots \\ b_{d}^{(1)} \end{pmatrix} \begin{pmatrix}b_{1}^{(2)}\\b_{2}^{(2)}\\ \vdots \\ b_{d}^{(2)} \end{pmatrix}\cdots \begin{pmatrix}b_{1}^{(p)}\\b_{2}^{(p)}\\ \vdots \\ b_{d}^{(p)} \end{pmatrix}.
$$
determines a $\frac{1}{k^{p}}$ hyper-cube $[0,\frac{1}{k^{p}}]^{d}+\{((0.b_{1}^{(1)}\cdots b_{1}^{(p)})_{k},\ldots, (0.b_{d}^{(1)}\cdots b_{d}^{(p)})_{k}))\}$. Then $\textbf{w}$ is in $\mathcal{L}_{X}$ if and only if this hyper-cube has non-empty intersection with some $X_{i}$\footnote{This formulation makes it evident that the entropy of $X$ will be equal to the box dimension of $\textbf{X}$. Since Hausdorff and box dimensions are often equal, this makes Theorem \ref{Theorem} quite plausible.}. Including both $k$-adic representations is required for this to be true. However, note that because we are tacitly assuming that the $k$-adic representations start $0.b_{1}\ldots$, we are only including one of the representations for each of $0$ and $1$, which is equivalent to excluding boxes which lie outside $[0,1]^{d}$.

We now count the number of level $p$ hyper-cubes which intersect each $X_{i}$. The number of level $1$ hyper-cubes intersecting $X_{i}$ is, by definition of the $A_{ij}$, equal to $\sum_{j}A_{ji}$. The number of level $2$-hyper-cubes intersecting $X_{i}$ is the sum of the number of times each $X_{j}$ occurs in $X_{i}$ times the sum of the number of times that each $X_{k}$ appears in those. That is $\sum_{k}\sum_{j}A_{kj}A_{ji}=\sum_{k}(A^{2})_{ki}$. Inductively, the number of level $p$ hyper-cubes intersecting $X_{i}$ is $\sum_{k}(A^{p})_{ki}= \mathbb{J}^{T} \cdot A^{p} \cdot e_{i}$, where $\mathbb{J}$ is the all ones vector. Each of these terms will be a lower bound for $\mathcal{P}(\mathcal{L}_{X},p)$. If we sum over $i$ this results in $\mathbb{J}^{T} \cdot A^{p} \cdot \mathbb{J}$. 
Because $A$ is non-negative, this is just the entrywise $1$-norm of $A^{p}$. By the above characterisations of $\mathcal{L}_{X}$, this quantity is an upper bound for $\mathcal{P}(\mathcal{L}_{X},p)$. Since there are $n$ terms in this sum, for each $p$ there must be at least one $i$ such that $\mathbb{J}^{T} \cdot A \cdot e_{i} \geqslant \frac{1}{n}\norm{A^{p}}$. By Gelfand's formula for the spectral radius, $\norm{A^{p}}_{1}=\rho(A)^{p}(1+o(1))^{p}$. This working has proven the following proposition.

\begin{proposition}\label{2}
Let $X$ be a $k$-self-similar set and $A$ the matrix $A=\{A_{ij}\}_{i,j=0}^{n-1}$. Then the entropy, $h(X)$, of $X$ is equal to
$$
h(X)=\lim_{p \to \infty}\frac{1}{p}\log_{k}(\mathcal{P}(\mathcal{L}_{X},p))=\log_{k}(\rho(A)).
$$
\end{proposition}

A comparison of Propositions \ref{1} and \ref{2} immediately yields Theorem \ref{Theorem}.

\section{An Extended Example}

In the hopes of helping the reader to parse the material presented here, we now illustrate the different parts of Theorem \ref{Theorem} in the context of a familiar example: the Cantor set, $\mathcal{C}$. This set can be constructed by starting with $[0,1]$ and repeatedly removing the open middle thirds from all remaining intervals.

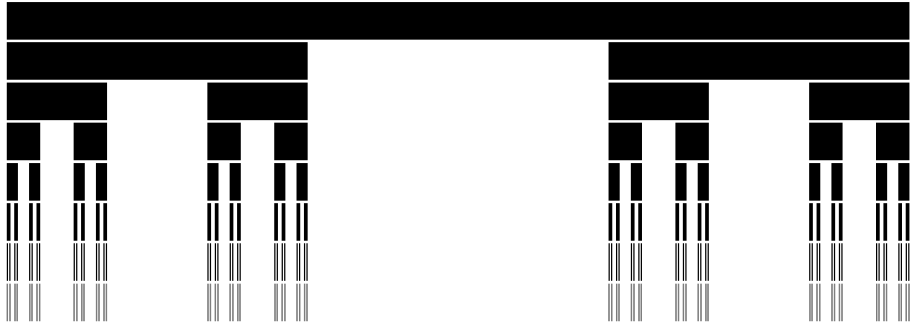
\begin{figure}[h]
    \centering
    \begin{tikzpicture}
    \CANTOR{0}{12}{0}
    \end{tikzpicture}
    \begin{tikzpicture}
    \CANTOR{1}{12}{0}
    \end{tikzpicture}
    \begin{tikzpicture}
    \CANTOR{2}{12}{0}
    \end{tikzpicture}
    \begin{tikzpicture}
    \CANTOR{3}{12}{0}
    \end{tikzpicture}
    \begin{tikzpicture}
    \CANTOR{4}{12}{0}
    \end{tikzpicture}
    \begin{tikzpicture}
    \CANTOR{5}{12}{0}
    \end{tikzpicture}
    \begin{tikzpicture}
    \CANTOR{6}{12}{0}
    \end{tikzpicture}
    \begin{tikzpicture}
    \CANTOR{7}{12}{0}
    \end{tikzpicture}
\end{figure}

The Cantor set is $3$-self-similar. Let us find its $3$-kernel. At the zeroth level, we intersect with $[0,1]$ and simply get $\mathcal{C}$. 

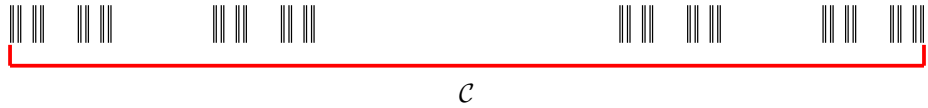
\begin{figure}[h!]
    \centering
\begin{tikzpicture}[yscale=0.555,xscale=0.45]

\CANTOR{6}{27}{0}

\draw[red,line width=0.05cm] (0,-1.0)--(0,-0.5);

\draw[red,line width=0.05cm] (27,-1.0)--(27,-0.5);

\draw[red,line width=0.05cm] (0,-1.0)--(27,-1.0);

\node(A) at (13.5,-1.0cm) {};
\node [below] at (A.south) {$\mathcal{C}$};

\end{tikzpicture}
\end{figure}

In the first level, we intersect with $[0,\frac{1}{3}],[\frac{1}{3},\frac{2}{3}]$ and $[\frac{2}{3},1]$. The first and last give $\mathcal{C}$, the middle $\{0,1\}$. 

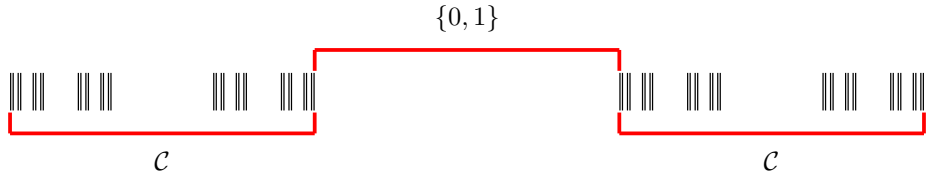
\begin{figure}[h!]
    \centering
\begin{tikzpicture}[yscale=0.555,xscale=0.45]

\CANTOR{6}{27}{0}

\draw[red,line width=0.05cm] (0,-1.0)--(0,-0.5);

\draw[red,line width=0.05cm] (9,-1.0)--(9,-0.5);

\draw[red,line width=0.05cm] (9,1.0)--(9,0.5);

\draw[red,line width=0.05cm] (18,1.0)--(18,0.5);

\draw[red,line width=0.05cm] (18,-1.0)--(18,-0.5);

\draw[red,line width=0.05cm] (27,-1.0)--(27,-0.5);

\draw[red,line width=0.05cm] (0,-1.0)--(9,-1.0);

\draw[red,line width=0.05cm] (9,1.0)--(18,1.0);

\draw[red,line width=0.05cm] (18,-1.0)--(27,-1.0);

\node(A) at (4.5,-1.0cm) {};
\node [below] at (A.south) {$\mathcal{C}$};

\node(B) at (13.5,1.0cm) {};
\node [above] at (B.north) {$\{0,1\}$};

\node(C) at (22.5,-1.0cm) {};
\node [below] at (C.south) {$\mathcal{C}$};

\end{tikzpicture}
\end{figure}

If we go to the second level, the interval in the middle gives $\{0\}$ and $\{1\}$.

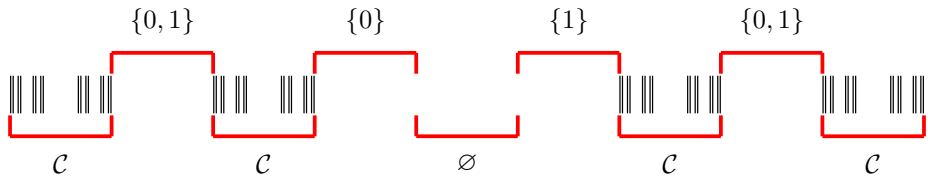
\begin{figure}[h!]
    \centering
\begin{tikzpicture}[yscale=0.555,xscale=0.45]

\CANTOR{6}{27}{0}

\draw[red,line width=0.05cm] (0,-1.0)--(0,-0.5);

\draw[red,line width=0.05cm] (3,-1.0)--(3,-0.5);

\draw[red,line width=0.05cm] (3,1.0)--(3,0.5);

\draw[red,line width=0.05cm] (6,1.0)--(6,0.5);

\draw[red,line width=0.05cm] (6,-1.0)--(6,-0.5);

\draw[red,line width=0.05cm] (9,-1.0)--(9,-0.5);

\draw[red,line width=0.05cm] (9,1.0)--(9,0.5);

\draw[red,line width=0.05cm] (12,1.0)--(12,0.5);

\draw[red,line width=0.05cm] (12,-1.0)--(12,-0.5);

\draw[red,line width=0.05cm] (15,-1.0)--(15,-0.5);

\draw[red,line width=0.05cm] (15,1.0)--(15,0.5);

\draw[red,line width=0.05cm] (18,1.0)--(18,0.5);

\draw[red,line width=0.05cm] (18,-1.0)--(18,-0.5);

\draw[red,line width=0.05cm] (21,-1.0)--(21,-0.5);

\draw[red,line width=0.05cm] (21,1.0)--(21,0.5);

\draw[red,line width=0.05cm] (24,1.0)--(24,0.5);

\draw[red,line width=0.05cm] (24,-1.0)--(24,-0.5);

\draw[red,line width=0.05cm] (27,-1.0)--(27,-0.5);

\draw[red,line width=0.05cm] (0,-1.0)--(3,-1.0);

\draw[red,line width=0.05cm] (3,1.0)--(6,1.0);

\draw[red,line width=0.05cm] (6,-1.0)--(9,-1.0);

\draw[red,line width=0.05cm] (9,1.0)--(12,1.0);

\draw[red,line width=0.05cm] (12,-1.0)--(15,-1.0);

\draw[red,line width=0.05cm] (15,1.0)--(18,1.0);

\draw[red,line width=0.05cm] (18,-1.0)--(21,-1.0);

\draw[red,line width=0.05cm] (21,1.0)--(24,1.0);

\draw[red,line width=0.05cm] (24,-1.0)--(27,-1.0);

\node(A) at (1.5,-1.0cm) {};
\node [below] at (A.south) {$\mathcal{C}$};

\node(B) at (4.5,1.0cm) {};
\node [above] at (B.north) {$\{0,1\}$};

\node(C) at (7.5,-1.0cm) {};
\node [below] at (C.south) {$\mathcal{C}$};

\node(D) at (10.5,1.0cm) {};
\node [above] at (D.north) {$\{0\}$};

\node(E) at (13.5,-1.0cm) {};
\node [below] at (E.south) {$\varnothing$};

\node(F) at (16.5,1.0cm) {};
\node [above] at (F.north) {$\{1\}$};

\node(G) at (19.5,-1.0cm) {};
\node [below] at (G.south) {$\mathcal{C}$};

\node(H) at (22.5,1.0cm) {};
\node [above] at (H.north) {$\{0,1\}$};

\node(I) at (25.5,-1.0cm) {};
\node [below] at (I.south) {$\mathcal{C}$};

\end{tikzpicture}
\end{figure}

From then on we get nothing new, and so we have $\mathcal{K}_{\mathcal{C}}^{3}=\{\mathcal{C},\{0,1\},\{0\},\{1\}\}$. The corresponding set $\textbf{X}=X_{0}\cup (X_{1}+\{2\})\cup (X_{2}+\{4\})\cup (X_{3}+\{6\})$ is in our case given by $\mathcal{C}\cup \{2,3,4,7\}$.

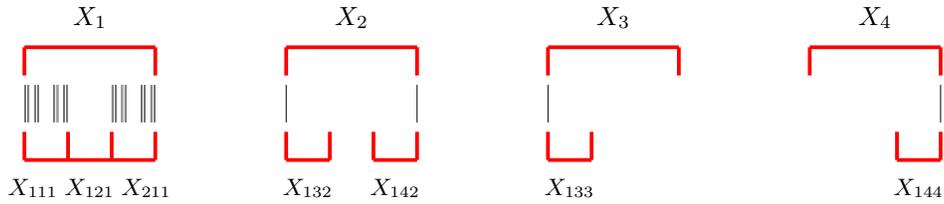
\begin{figure}[h!]
    \centering
\begin{tikzpicture}[yscale=0.75,xscale=0.58]

\CANTOR{5}{3}{0}

\draw[line width=0.01cm] (6,-0.333)--(6,0.333);
\draw[line width=0.01cm] (9,-0.333)--(9,0.333);
\draw[line width=0.01cm] (12,-0.333)--(12,0.333);
\draw[line width=0.01cm] (21,-0.333)--(21,0.333);

\draw[red,line width=0.05cm] (0,1.0)--(0,0.5);
\draw[red,line width=0.05cm] (3,1.0)--(3,0.5);
\draw[red,line width=0.05cm] (0,1.0)--(3,1.0);

\draw[red,line width=0.05cm] (6,1.0)--(6,0.5);
\draw[red,line width=0.05cm] (9,1.0)--(9,0.5);
\draw[red,line width=0.05cm] (6,1.0)--(9,1.0);

\draw[red,line width=0.05cm] (12,1.0)--(12,0.5);
\draw[red,line width=0.05cm] (15,1.0)--(15,0.5);
\draw[red,line width=0.05cm] (12,1.0)--(15,1.0);

\draw[red,line width=0.05cm] (18,1.0)--(18,0.5);
\draw[red,line width=0.05cm] (21,1.0)--(21,0.5);
\draw[red,line width=0.05cm] (18,1.0)--(21,1.0);

\node(A) at (1.5,1.0cm) {};
\node [above] at (A.north) {$X_{1}$};
\node(B) at (7.5,1.0cm) {};
\node [above] at (B.north) {$X_{2}$};
\node(C) at (13.5,1.0cm) {};
\node [above] at (C.north) {$X_{3}$};
\node(D) at (19.5,1.0cm) {};
\node [above] at (D.north) {$X_{4}$};

\draw[red,line width=0.05cm] (0,-1.0)--(0,-0.5);
\draw[red,line width=0.05cm] (1,-1.0)--(1,-0.5);
\draw[red,line width=0.05cm] (2,-1.0)--(2,-0.5);
\draw[red,line width=0.05cm] (3,-1.0)--(3,-0.5);
\draw[red,line width=0.05cm] (0,-1.0)--(3,-1.0);

\node(E) at (0.2,-1.0cm) {};
\node [below] at (E.south) {\small $X_{111}$};
\node(F) at (1.5,-1.0cm) {};
\node [below] at (F.south) {\small $X_{121}$};
\node(G) at (2.8,-1.0cm) {};
\node [below] at (G.south) {\small $X_{211}$};

\draw[red,line width=0.05cm] (6,-1.0)--(6,-0.5);
\draw[red,line width=0.05cm] (7,-1.0)--(7,-0.5);
\draw[red,line width=0.05cm] (8,-1.0)--(8,-0.5);
\draw[red,line width=0.05cm] (9,-1.0)--(9,-0.5);
\draw[red,line width=0.05cm] (6,-1.0)--(7,-1.0);
\draw[red,line width=0.05cm] (8,-1.0)--(9,-1.0);

\node(H) at (6.5,-1.0cm) {};
\node [below] at (H.south) {\small $X_{132}$};
\node(I) at (8.5,-1.0cm) {};
\node [below] at (I.south) {\small $X_{142}$};

\draw[red,line width=0.05cm] (12,-1.0)--(12,-0.5);
\draw[red,line width=0.05cm] (13,-1.0)--(13,-0.5);
\draw[red,line width=0.05cm] (12,-1.0)--(13,-1.0);

\node(J) at (12.5,-1.0cm) {};
\node [below] at (J.south) {\small $X_{133}$};

\draw[red,line width=0.05cm] (20,-1.0)--(20,-0.5);
\draw[red,line width=0.05cm] (21,-1.0)--(21,-0.5);
\draw[red,line width=0.05cm] (20,-1.0)--(21,-1.0);

\node(K) at (20.5,-1.0cm) {};
\node [below] at (K.south) {\small $X_{144}$};

\end{tikzpicture}
\end{figure}

Examining this, we see that $\mathcal{C}$ contains two $\frac{1}{3}$ scale copies of $\mathcal{C}$ and one copy of $\{0,1\}$, and that $\{0,1\}$ contains one copy each of $\{0\}$ and $\{1\}$ which only contain copies of themselves. Hence the $A$ matrix is
$$
A=\begin{bmatrix} 2 & 0 & 0 & 0 \\ 1 & 0 & 0 & 0\\ 0 & 1 & 1 & 0\\ 0 & 1 & 0 & 1 \\ \end{bmatrix}.
$$
Since this is lower triangular, the diagonal entries are its eigenvalues, so the spectral radius is $2$. The dimension of the Cantor set is thus $\log_{3}(2)$, a well known result.

The next diagram shows the graph $G$. Recall that the vertices are essentially the $X_{hij}$, the $\frac{1}{3}$ scaled copies of the $X_{i}$ contained in $\textbf{X}$, and each small copy of $X_{i}$ is connected to every element of the full scale copy of $X_{i}$.
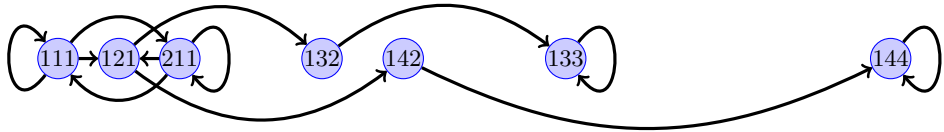
\begin{figure}[h]
    \centering
\begin{tikzpicture}[yscale=0.75,xscale=0.54]
\node(E) at (0,-1.0cm) [circle,draw=blue,fill=blue!20,inner sep=0.1] {\small $111$};
\node(F) at (1.5,-1.0cm) [circle,draw=blue,fill=blue!20,inner sep=0.1] {\small $121$};
\node(G) at (3,-1.0cm) [circle,draw=blue,fill=blue!20,inner sep=0.1] {\small $211$};
\node(H) at (6.5,-1.0cm) [circle,draw=blue,fill=blue!20,inner sep=0.1] {\small $132$};
\node(I) at (8.5,-1.0cm) [circle,draw=blue,fill=blue!20,inner sep=0.1] {\small $142$};
\node(J) at (12.5,-1.0cm) [circle,draw=blue,fill=blue!20,inner sep=0.1] {\small $133$};
\node(K) at (20.5,-1.0cm) [circle,draw=blue,fill=blue!20,inner sep=0.1] {\small $144$};

\draw [->,very thick] (E) to (F);
\draw [->,very thick] (G) to (F);
\draw [->,very thick] (E) to [out=45,in=135] (G);
\draw [->,very thick] (G) to [out=-135,in=-45] (E);
\draw [->,very thick] (E) to [out=-135,in=135, looseness=8] (E);
\draw [->,very thick] (G) to [out=45,in=-45, looseness=8] (G);
\draw [->,very thick] (F) to [out=35,in=145] (H);
\draw [->,very thick] (F) to [out=-30,in=-150] (I);
\draw [->,very thick] (H) to [out=30,in=150] (J);
\draw [->,very thick] (I) to [out=-20,in=-160] (K);
\draw [->,very thick] (J) to [out=45,in=-45, looseness=8] (J);
\draw [->,very thick] (K) to [out=45,in=-45, looseness=8] (K);

\end{tikzpicture}
\end{figure}

The corresponding similarity maps are demonstrated in the next diagram. Recall that $T_{(hij),(rsi)}$ is given by the translation and scaling which maps the large interval containing $X_{i}$ bijectively to the length $\frac{1}{3}$ interval which contains $X_{hij}$. This is independent of $(rsi)$, so we let $T_{(hij),(rsi)}=T_{hij}$.

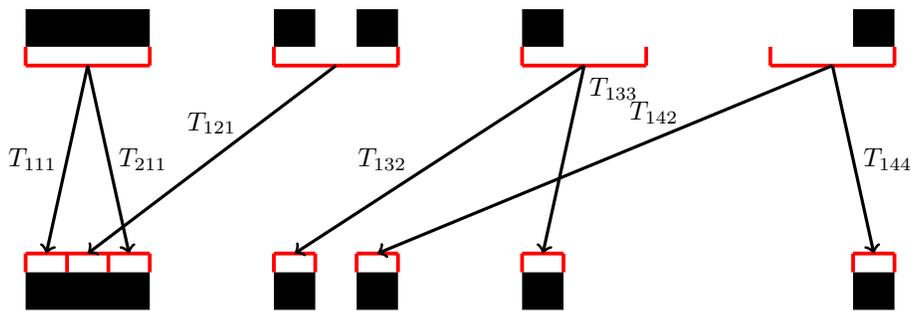
\begin{figure}[h!]
    \centering
\begin{tikzpicture}[xscale=0.55, yscale=0.5]

\begin{scope}
\draw[line width=0.5cm] (0,0)--(3,0);
\draw[line width=0.5cm] (6,0)--(7,0);
\draw[line width=0.5cm] (8,0)--(9,0);
\draw[line width=0.5cm] (12,0)--(13,0);
\draw[line width=0.5cm] (20,0)--(21,0);

\draw[red,line width=0.05cm] (0,-1.0)--(0,-0.5);
\draw[red,line width=0.05cm] (3,-1.0)--(3,-0.5);
\draw[red,line width=0.05cm] (0,-1.0)--(3,-1.0);

\draw[red,line width=0.05cm] (6,-1.0)--(6,-0.5);
\draw[red,line width=0.05cm] (9,-1.0)--(9,-0.5);
\draw[red,line width=0.05cm] (6,-1.0)--(9,-1.0);

\draw[red,line width=0.05cm] (12,-1.0)--(12,-0.5);
\draw[red,line width=0.05cm] (15,-1.0)--(15,-0.5);
\draw[red,line width=0.05cm] (12,-1.0)--(15,-1.0);

\draw[red,line width=0.05cm] (18,-1.0)--(18,-0.5);
\draw[red,line width=0.05cm] (21,-1.0)--(21,-0.5);
\draw[red,line width=0.05cm] (18,-1.0)--(21,-1.0);

\node(A) at (1.5,-1.0cm) {};
\node(B) at (7.5,-1.0cm) {};
\node(C) at (13.5,-1.0cm) {};
\node(D) at (19.5,-1.0cm) {};

\draw[red,line width=0.05cm] (0,-6.0)--(0,-6.5);
\draw[red,line width=0.05cm] (1,-6.0)--(1,-6.5);
\draw[red,line width=0.05cm] (2,-6.0)--(2,-6.5);
\draw[red,line width=0.05cm] (3,-6.0)--(3,-6.5);
\draw[red,line width=0.05cm] (0,-6.0)--(3,-6.0);
\draw[red,line width=0.05cm] (6,-6.0)--(6,-6.5);
\draw[red,line width=0.05cm] (7,-6.0)--(7,-6.5);
\draw[red,line width=0.05cm] (8,-6.0)--(8,-6.5);
\draw[red,line width=0.05cm] (9,-6.0)--(9,-6.5);
\draw[red,line width=0.05cm] (6,-6.0)--(7,-6.0);
\draw[red,line width=0.05cm] (8,-6.0)--(9,-6.0);
\draw[red,line width=0.05cm] (12,-6.0)--(12,-6.5);
\draw[red,line width=0.05cm] (13,-6.0)--(13,-6.5);
\draw[red,line width=0.05cm] (12,-6.0)--(13,-6.0);
\draw[red,line width=0.05cm] (20,-6.0)--(20,-6.5);
\draw[red,line width=0.05cm] (21,-6.0)--(21,-6.5);
\draw[red,line width=0.05cm] (20,-6.0)--(21,-6.0);

\node(E) at (0.5,-6.0cm) {};
\node(F) at (1.5,-6.0cm) {};
\node(G) at (2.5,-6.0cm) {};
\node(H) at (6.5,-6.0cm) {};
\node(I) at (8.5,-6.0cm) {};
\node(J) at (12.5,-6.0cm) {};
\node(K) at (20.5,-6.0cm) {};

\draw [->,very thick] (A.center) -- node[left] {$T_{111}$} (E.center);
\draw [->,very thick] (A.center) -- node[right] {$T_{211}$} (G.center);
\draw [->,very thick] (B.center) -- node[above=0.2cm] {$T_{121}$} (F.center);
\draw [->,very thick] (C.center) -- node[left=0.3cm] {$T_{132}$} (H.center);
\draw [->,very thick] (C.center) -- node[very near start, right] {$T_{133}$} (J.center);
\draw [->,very thick] (D.center) -- node[near start,left=0.4cm] {$T_{142}$} (I.center);
\draw [->,very thick] (D.center) -- node[right] {$T_{144}$}(K.center);
\end{scope}

\begin{scope}[yshift=-7.0cm]

\draw[line width=0.5cm] (0,0)--(3,0);
\draw[line width=0.5cm] (6,0)--(7,0);
\draw[line width=0.5cm] (8,0)--(9,0);
\draw[line width=0.5cm] (12,0)--(13,0);
\draw[line width=0.5cm] (20,0)--(21,0);
\end{scope}
\end{tikzpicture}
\end{figure}

This construction is not very efficient; evidently only the parts not relating to $\{0,1\},\{0\}$ and $\{1\}$ are needed to construct $\mathcal{C}$. Our method will often have these vestigial pieces relating to the `boundary' of the set which are not necessary if one is only interested in constructing $X$. However, the more complicated setup is guaranteed to work and simplifies the book-keeping required in the proof. 

The next diagram shows the zeroth, first, second and third levels of the construction, showing what maps where between each stage. This may help in visualising which $p$-digit $3$-adic rationals have their corresponding $\frac{1}{3^{p}}$ interval intersecting $\textbf{X}$, which helps in understanding the following results about $\mathcal{L}_{\mathcal{C}}$.

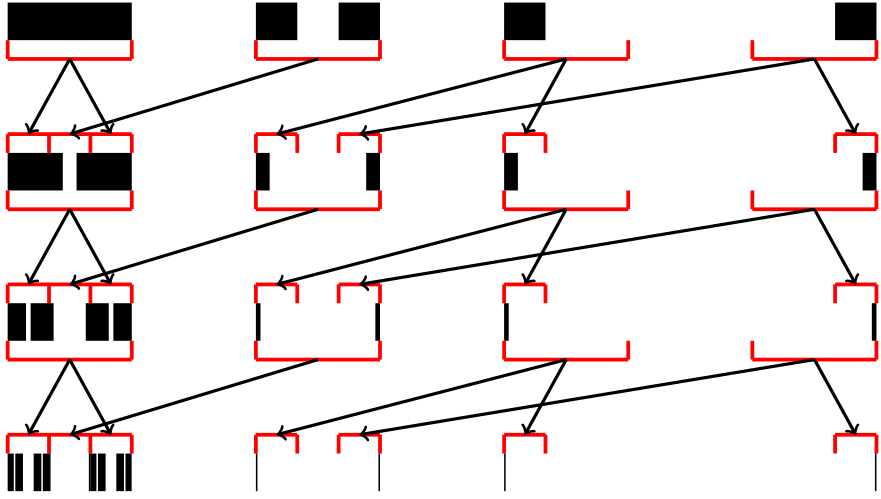
\begin{figure}[h]
    \centering
\begin{tikzpicture}[xscale=0.55,yscale=0.5]

\begin{scope}
\draw[line width=0.5cm] (0,0)--(3,0);
\draw[line width=0.5cm] (6,0)--(7,0);
\draw[line width=0.5cm] (8,0)--(9,0);
\draw[line width=0.5cm] (12,0)--(13,0);
\draw[line width=0.5cm] (20,0)--(21,0);

\draw[red,line width=0.05cm] (0,-1.0)--(0,-0.5);
\draw[red,line width=0.05cm] (3,-1.0)--(3,-0.5);
\draw[red,line width=0.05cm] (0,-1.0)--(3,-1.0);

\draw[red,line width=0.05cm] (6,-1.0)--(6,-0.5);
\draw[red,line width=0.05cm] (9,-1.0)--(9,-0.5);
\draw[red,line width=0.05cm] (6,-1.0)--(9,-1.0);

\draw[red,line width=0.05cm] (12,-1.0)--(12,-0.5);
\draw[red,line width=0.05cm] (15,-1.0)--(15,-0.5);
\draw[red,line width=0.05cm] (12,-1.0)--(15,-1.0);

\draw[red,line width=0.05cm] (18,-1.0)--(18,-0.5);
\draw[red,line width=0.05cm] (21,-1.0)--(21,-0.5);
\draw[red,line width=0.05cm] (18,-1.0)--(21,-1.0);

\node(A) at (1.5,-1.0cm) {};
\node(B) at (7.5,-1.0cm) {};
\node(C) at (13.5,-1.0cm) {};
\node(D) at (19.5,-1.0cm) {};

\draw[red,line width=0.05cm] (0,-3.0)--(0,-3.5);
\draw[red,line width=0.05cm] (1,-3.0)--(1,-3.5);
\draw[red,line width=0.05cm] (2,-3.0)--(2,-3.5);
\draw[red,line width=0.05cm] (3,-3.0)--(3,-3.5);
\draw[red,line width=0.05cm] (0,-3.0)--(3,-3.0);
\draw[red,line width=0.05cm] (6,-3.0)--(6,-3.5);
\draw[red,line width=0.05cm] (7,-3.0)--(7,-3.5);
\draw[red,line width=0.05cm] (8,-3.0)--(8,-3.5);
\draw[red,line width=0.05cm] (9,-3.0)--(9,-3.5);
\draw[red,line width=0.05cm] (6,-3.0)--(7,-3.0);
\draw[red,line width=0.05cm] (8,-3.0)--(9,-3.0);
\draw[red,line width=0.05cm] (12,-3.0)--(12,-3.5);
\draw[red,line width=0.05cm] (13,-3.0)--(13,-3.5);
\draw[red,line width=0.05cm] (12,-3.0)--(13,-3.0);
\draw[red,line width=0.05cm] (20,-3.0)--(20,-3.5);
\draw[red,line width=0.05cm] (21,-3.0)--(21,-3.5);
\draw[red,line width=0.05cm] (20,-3.0)--(21,-3.0);

\node(E) at (0.5,-3.0cm) {};
\node(F) at (1.5,-3.0cm) {};
\node(G) at (2.5,-3.0cm) {};
\node(H) at (6.5,-3.0cm) {};
\node(I) at (8.5,-3.0cm) {};
\node(J) at (12.5,-3.0cm) {};
\node(K) at (20.5,-3.0cm) {};

\draw [->,very thick] (A.center) to (E.center);
\draw [->,very thick] (A.center) to (G.center);
\draw [->,very thick] (B.center) to (F.center);
\draw [->,very thick] (C.center) to (H.center);
\draw [->,very thick] (C.center) to (J.center);
\draw [->,very thick] (D.center) to (I.center);
\draw [->,very thick] (D.center) to (K.center);
\end{scope}

\begin{scope}[yshift=-4.0cm]

\draw[line width=0.5cm] (0,0)--(1.333,0);
\draw[line width=0.5cm] (1.666,0)--(3,0);
\draw[line width=0.5cm] (6,0)--(6.333,0);
\draw[line width=0.5cm] (8.666,0)--(9,0);
\draw[line width=0.5cm] (12,0)--(12.333,0);
\draw[line width=0.5cm] (20.666,0)--(21,0);

\draw[red,line width=0.05cm] (0,-1.0)--(0,-0.5);
\draw[red,line width=0.05cm] (3,-1.0)--(3,-0.5);
\draw[red,line width=0.05cm] (0,-1.0)--(3,-1.0);

\draw[red,line width=0.05cm] (6,-1.0)--(6,-0.5);
\draw[red,line width=0.05cm] (9,-1.0)--(9,-0.5);
\draw[red,line width=0.05cm] (6,-1.0)--(9,-1.0);

\draw[red,line width=0.05cm] (12,-1.0)--(12,-0.5);
\draw[red,line width=0.05cm] (15,-1.0)--(15,-0.5);
\draw[red,line width=0.05cm] (12,-1.0)--(15,-1.0);

\draw[red,line width=0.05cm] (18,-1.0)--(18,-0.5);
\draw[red,line width=0.05cm] (21,-1.0)--(21,-0.5);
\draw[red,line width=0.05cm] (18,-1.0)--(21,-1.0);

\node(A) at (1.5,-1.0cm) {};
\node(B) at (7.5,-1.0cm) {};
\node(C) at (13.5,-1.0cm) {};
\node(D) at (19.5,-1.0cm) {};

\draw[red,line width=0.05cm] (0,-3.0)--(0,-3.5);
\draw[red,line width=0.05cm] (1,-3.0)--(1,-3.5);
\draw[red,line width=0.05cm] (2,-3.0)--(2,-3.5);
\draw[red,line width=0.05cm] (3,-3.0)--(3,-3.5);
\draw[red,line width=0.05cm] (0,-3.0)--(3,-3.0);
\draw[red,line width=0.05cm] (6,-3.0)--(6,-3.5);
\draw[red,line width=0.05cm] (7,-3.0)--(7,-3.5);
\draw[red,line width=0.05cm] (8,-3.0)--(8,-3.5);
\draw[red,line width=0.05cm] (9,-3.0)--(9,-3.5);
\draw[red,line width=0.05cm] (6,-3.0)--(7,-3.0);
\draw[red,line width=0.05cm] (8,-3.0)--(9,-3.0);
\draw[red,line width=0.05cm] (12,-3.0)--(12,-3.5);
\draw[red,line width=0.05cm] (13,-3.0)--(13,-3.5);
\draw[red,line width=0.05cm] (12,-3.0)--(13,-3.0);
\draw[red,line width=0.05cm] (20,-3.0)--(20,-3.5);
\draw[red,line width=0.05cm] (21,-3.0)--(21,-3.5);
\draw[red,line width=0.05cm] (20,-3.0)--(21,-3.0);

\node(E) at (0.5,-3.0cm) {};
\node(F) at (1.5,-3.0cm) {};
\node(G) at (2.5,-3.0cm) {};
\node(H) at (6.5,-3.0cm) {};
\node(I) at (8.5,-3.0cm) {};
\node(J) at (12.5,-3.0cm) {};
\node(K) at (20.5,-3.0cm) {};

\draw [->,very thick] (A.center) to (E.center);
\draw [->,very thick] (A.center) to (G.center);
\draw [->,very thick] (B.center) to (F.center);
\draw [->,very thick] (C.center) to (H.center);
\draw [->,very thick] (C.center) to (J.center);
\draw [->,very thick] (D.center) to (I.center);
\draw [->,very thick] (D.center) to (K.center);

\end{scope}

\begin{scope}[yshift=-8.0cm]

\draw[line width=0.5cm] (0,0)--(0.444,0);
\draw[line width=0.5cm] (0.555,0)--(1.111,0);
\draw[line width=0.5cm] (1.888,0)--(2.444,0);
\draw[line width=0.5cm] (2.555,0)--(3,0);

\draw[line width=0.5cm] (6,0)--(6.111,0);
\draw[line width=0.5cm] (8.888,0)--(9,0);
\draw[line width=0.5cm] (12,0)--(12.111,0);
\draw[line width=0.5cm] (20.888,0)--(21,0);

\draw[red,line width=0.05cm] (0,-1.0)--(0,-0.5);
\draw[red,line width=0.05cm] (3,-1.0)--(3,-0.5);
\draw[red,line width=0.05cm] (0,-1.0)--(3,-1.0);

\draw[red,line width=0.05cm] (6,-1.0)--(6,-0.5);
\draw[red,line width=0.05cm] (9,-1.0)--(9,-0.5);
\draw[red,line width=0.05cm] (6,-1.0)--(9,-1.0);

\draw[red,line width=0.05cm] (12,-1.0)--(12,-0.5);
\draw[red,line width=0.05cm] (15,-1.0)--(15,-0.5);
\draw[red,line width=0.05cm] (12,-1.0)--(15,-1.0);

\draw[red,line width=0.05cm] (18,-1.0)--(18,-0.5);
\draw[red,line width=0.05cm] (21,-1.0)--(21,-0.5);
\draw[red,line width=0.05cm] (18,-1.0)--(21,-1.0);

\node(A) at (1.5,-1.0cm) {};
\node(B) at (7.5,-1.0cm) {};
\node(C) at (13.5,-1.0cm) {};
\node(D) at (19.5,-1.0cm) {};

\draw[red,line width=0.05cm] (0,-3.0)--(0,-3.5);
\draw[red,line width=0.05cm] (1,-3.0)--(1,-3.5);
\draw[red,line width=0.05cm] (2,-3.0)--(2,-3.5);
\draw[red,line width=0.05cm] (3,-3.0)--(3,-3.5);
\draw[red,line width=0.05cm] (0,-3.0)--(3,-3.0);
\draw[red,line width=0.05cm] (6,-3.0)--(6,-3.5);
\draw[red,line width=0.05cm] (7,-3.0)--(7,-3.5);
\draw[red,line width=0.05cm] (8,-3.0)--(8,-3.5);
\draw[red,line width=0.05cm] (9,-3.0)--(9,-3.5);
\draw[red,line width=0.05cm] (6,-3.0)--(7,-3.0);
\draw[red,line width=0.05cm] (8,-3.0)--(9,-3.0);
\draw[red,line width=0.05cm] (12,-3.0)--(12,-3.5);
\draw[red,line width=0.05cm] (13,-3.0)--(13,-3.5);
\draw[red,line width=0.05cm] (12,-3.0)--(13,-3.0);
\draw[red,line width=0.05cm] (20,-3.0)--(20,-3.5);
\draw[red,line width=0.05cm] (21,-3.0)--(21,-3.5);
\draw[red,line width=0.05cm] (20,-3.0)--(21,-3.0);

\node(E) at (0.5,-3.0cm) {};
\node(F) at (1.5,-3.0cm) {};
\node(G) at (2.5,-3.0cm) {};
\node(H) at (6.5,-3.0cm) {};
\node(I) at (8.5,-3.0cm) {};
\node(J) at (12.5,-3.0cm) {};
\node(K) at (20.5,-3.0cm) {};

\draw [->,very thick] (A.center) to (E.center);
\draw [->,very thick] (A.center) to (G.center);
\draw [->,very thick] (B.center) to (F.center);
\draw [->,very thick] (C.center) to (H.center);
\draw [->,very thick] (C.center) to (J.center);
\draw [->,very thick] (D.center) to (I.center);
\draw [->,very thick] (D.center) to (K.center);

\end{scope}

\begin{scope}[yshift=-12.0cm]

\begin{scope}[xscale=0.333]
\draw[line width=0.5cm] (0,0)--(0.444,0);
\draw[line width=0.5cm] (0.555,0)--(1.111,0);
\draw[line width=0.5cm] (1.888,0)--(2.444,0);
\draw[line width=0.5cm] (2.555,0)--(3,0);
\end{scope}

\begin{scope}[xscale=0.333,xshift=6.0cm]
\draw[line width=0.5cm] (0,0)--(0.444,0);
\draw[line width=0.5cm] (0.555,0)--(1.111,0);
\draw[line width=0.5cm] (1.888,0)--(2.444,0);
\draw[line width=0.5cm] (2.555,0)--(3,0);
\end{scope}

\draw[line width=0.5cm] (1,0)--(1.037,0);
\draw[line width=0.5cm] (1.963,0)--(2,0);

\draw[line width=0.5cm] (6,0)--(6.037,0);
\draw[line width=0.5cm] (8.963,0)--(9,0);
\draw[line width=0.5cm] (12,0)--(12.037,0);
\draw[line width=0.5cm] (20.963,0)--(21,0);

\end{scope}

\end{tikzpicture}
\end{figure}

The language, $\mathcal{L}_{\mathcal{C}}$ consists of all finite words which appear in the $3$-adic expansion of a real number in $\mathcal{C}$ (and which start $0.b_{1}\ldots$). It is well known that $\mathcal{C}$ contains all real numbers in $[0,1]$ which have a base $3$ representation which contains only $0$'s and $2$'s. Hence $\widetilde{\mathcal{L}}=\{a_{1}a_{2}\cdots a_{p}|a_{i}=0,2\}$ is contained in $\mathcal{L}_{\mathcal{C}}$. There are exactly $2^{p}$ words of length $p$ in here, and so the entropy is easily calculated to be $\log_{3}(2)$. However, $\widetilde{\mathcal{L}}$ is not all of $\mathcal{L}_{\mathcal{C}}$. We must also include the words corresponding to the double representations of endpoints of the intervals of the set. These correspond to words which begin on the alphabet $0,2$, then have a $1$, and then have trailing $2$'s or $0$'s. We count these by starting with a length $p-k$ $\widetilde{\mathcal{L}}$ word, then appending a $1$ and then $k-1$ $0$'s or $k-1$ $2$'s). This is an extra
$$
2^{p-1}+2(2^{p-2}+2^{p-3}+\cdots+2+1)=3\cdot 2^{p-1}-2
$$
words (for $p\geqslant 2$ at least). Thus the total number of length $p$ words in $\mathcal{L}_{\mathcal{C}}$ is $5 \cdot 2^{p-1}-2$, and the entropy is still $\log_{3}(2)$. This exhibits the general fact that including these multiple representations will not change the entropy because it only duplicates existing words, and hence increases the size of the language by at most a constant multiple, not changing the asymptotic growth of $\log_{3}(\mathcal{P}(\mathcal{L},p))$.

Interestingly enough, in this case is turns out that $\mathcal{P}(\mathcal{L}_{\mathcal{C}},p)=\mathbb{J}^{T} \cdot A^{p} \cdot e_{1}$. That is, all words of length $p$ which appear in $\mathcal{L}_{\mathcal{C}}$ appear as the first $p$ ternary digits of some expansion of some real number in $\mathcal{C}$. After $A=A^{1}$ there is an easily observable pattern to the $A^{p}$ matrices which makes this equality obvious. We leave it up to the reader to interpret what each of the terms in the first column of these matrices and what they say about where all of the different words in $\mathcal{L}_{\mathcal{C}}$ are coming from.

\small
$$
A=\begin{bmatrix} 2 & 0 & 0 & 0 \\ 1 & 0 & 0 & 0\\ 0 & 1 & 1 & 0\\ 0 & 1 & 0 & 1 \\ \end{bmatrix}\; A^{2}=\begin{bmatrix} 4 & 0 & 0 & 0 \\ 2 & 0 & 0 & 0\\ 1 & 1 & 1 & 0\\ 1 & 1 & 0 & 1 \\  \end{bmatrix} \; A^{3}=\begin{bmatrix} 8 & 0 & 0 & 0 \\ 4 & 0 & 0 & 0\\ 3 & 1 & 1 & 0\\ 3 & 1 & 0 & 1 \\  \end{bmatrix} \; A^{4}=\begin{bmatrix} 16 & 0 & 0 & 0 \\ 8 & 0 & 0 & 0\\ 7 & 1 & 1 & 0\\ 7 & 1 & 0 & 1 \\  \end{bmatrix} \cdots
$$

\section{A final comment}

Theorem \ref{Theorem} proves that $\dim_{\mathcal{H}}(X)=\log_{k}(\rho)$, where $\rho$ is some algebraic number. This quantity is some kind of scaling exponent, telling us how the complexity of the object grows as one tends towards some limiting place. A very similar quantity occurs in the study of Mahler functions and $k$-regular sequences. If $f(x)$ is $k$-Mahler and the associated series converges in the unit disk, then one can prove that
$$f(x) \asymp \frac{1}{(1-x)^{\log_{k}(\rho)}}$$
as $x \to 1$, where $\rho$ is some algebraic number \cite{BellCoonsMahlerEigenvalue1}. For a non-negative $k$-regular sequence $f(n)$, one can prove that

$$
\sum_{n \leqslant N} f(n) \asymp N^{\log_{k}(\rho)} \log(N)^{a}
$$

The presence of a $\log_{k}(\rho)$ type exponent is another common feature in this circle of related objects, in addition to the Cobham type independence results. One might wonder if this quantity can be connected to the independence results. If an object was both $k$- and $l$-regular/Mahler/self-similar for multiplicatively independent $k$ and $l$, then because this exponent is defined independently of $k$ and $l$, we must have an equality of the form $\log_{k}(\alpha)=\log_{l}(\beta)$ for some algebraic $\alpha$ and $\beta$. It is a well known heuristic in number theory that the logarithms of algebraic numbers to different bases are `independent' in some sense, and so equalities like $\log_{k}(\alpha)=\log_{l}(\beta)$ are unlikely to be possible.

Unfortunately, our current understanding of when equalities like this can occur is incomplete, so this idea cannot be used to prove the desired independence results results. Exactly when we can have
$$
\log(\alpha)\log(l)=\log(\beta)\log(k)
$$
is the content of the famous, but as yet unproven, four exponentials conjecture.

\bibliographystyle{amsplain}

\end{document}